\documentclass[a4paper,12pt,reqno]{amsart}

\usepackage[english]{babel}
\usepackage[T1]{fontenc}
\usepackage[left=1in,right=1in,top=1.2in,bottom=1.2in]{geometry}
\usepackage{times}
\usepackage{microtype}

\usepackage{commath,amssymb,amscd,mathrsfs,mathtools,dsfont,bbm,multirow,array,caption,comment,pifont}
\usepackage[all]{xy}
\usepackage{enumerate}
\usepackage{longtable}

\usepackage{setspace}
\allowdisplaybreaks

\usepackage[dvipsnames,svgnames,table]{xcolor}
\definecolor{red}{RGB}{255,25,25}
\definecolor{blue}{RGB}{25,50,200}
\usepackage{tikz-cd}
\usepackage[pagebackref,linktocpage]{hyperref}

\hypersetup{
pdftitle={},				
pdfauthor={Hu--Truong--Xie},	
colorlinks=true,			
linkcolor=red,			
citecolor=MidnightBlue,	
filecolor=magenta,		
urlcolor=MidnightBlue	
}

\usepackage{cleveref}	

\newtheorem{theorem}{Theorem}[section]
\crefname{theorem}{Theorem}{Theorems}
\newtheorem{lemma}[theorem]{Lemma}
\crefname{lemma}{Lemma}{Lemmas}
\newtheorem{proposition}[theorem]{Proposition}
\crefname{proposition}{Proposition}{Propositions}

\crefname{prop}{Proposition}{Propositions}

\crefname{corollary}{Corollary}{Corollaries}

\crefname{cor}{Corollary}{Corollaries}

\crefname{conjecture}{Conjecture}{Conjectures}

\crefname{conj}{Conjecture}{Conjectures}
\newtheorem*{conj*}{Conjecture}
\crefname{conj}{Conjecture}{Conjectures}

\crefname{conjA}{Conjecture}{Conjecture}

\crefname{conjB}{Conjecture}{Conjecture}

\crefname{conjC}{Conjecture}{Conjecture}

\crefname{conjDk}{Conjecture}{Conjecture}

\crefname{conjD}{Conjecture}{Conjecture}

\crefname{conjH}{Conjecture}{Conjecture}

\crefname{conjGr}{Conjecture}{Conjecture}

\theoremstyle{definition}

\crefname{definition}{Definition}{Definitions}

\crefname{defn}{Definition}{Definitions}

\crefname{example}{Example}{Examples}

\crefname{notation}{Notation}{Notation}
\newtheorem*{notation*}{Notation}
\crefname{notation}{Notation}{Notation}

\crefname{problem}{Problem}{Problems}

\crefname{question}{Question}{Questions}

\crefname{condition}{Condition}{Conditions}

\crefname{assumption}{Assumption}{Assumptions}

\crefname{propGr}{Property}{Property}

\theoremstyle{remark}

\crefname{rmk}{Remark}{Remarks}
\newtheorem*{rmk*}{Remark}
\crefname{rmk}{Remark}{Remarks}

\crefname{remark}{Remark}{Remarks}

\crefname{fact}{Fact}{Facts}

\crefname{claim}{Claim}{Claims}
\newtheorem*{claim*}{Claim}
\crefname{claim}{Claim}{Claims}

\crefname{step}{Step}{Steps}

\crefname{case}{Case}{Cases}

\numberwithin{equation}{section}

\newcommand{\doi}[1]{\href{https://doi.org/#1}{{\tt DOI:#1}}}

\newcommand{\bZ}{\mathbf{Z}}

\title[Standard conjecture D and some conjectures around Weil's Riemann hypothesis]{Standard conjecture D and some conjectures around Weil's Riemann hypothesis}
   \author{Tuyen Trung Truong}
 \address{Department of mathematics, University of Oslo, Norway}
 \email{tuyentt@math.uio.no}

\thanks{}
    \date{\today}
    \keywords{Correspondence, $l$-adic cohomology, Polarised endomorphism, Rational map, Relative dynamical degrees, Spectral radius, Standard conjectures, Semi-simplicity, Weil's Riemann hypothesis}
    \subjclass[2010]{37F, 14D, 32U40, 32H50}

\begin{document}
\maketitle
\begin{abstract}



Let $X$ be a smooth projective variety defined on a finite field $\mathbb{F}_q$. On $X$ there is a special morphism $Fr_X$, which raises coordinates to exponent $q$: $t\mapsto t^q$. The two main results in this paper are: 

Result 1: If Standard conjecture D holds (for algebraic cycles of dimension $=\dim (X)$) on $X\times X$, then all polarised endomorphisms on $X$ are semisimple. 

Result 2: We provide heuristic arguments to show that Standard Conjecture D should imply both Dynamical degree comparison conjecture (a generalisation of both Weil's Riemann hypothesis and Tate's question on the absolute value of the eigenvalues of polarised endomorphisms),  Norm comparison conjecture (allowing to bound the growth  of the pullback of iterations of an endomorphism on cohomology groups in terms of that on algebraic cycles, in particularly implying the semisimplicity of polarised endomorphisms), and Conjecture $G_r$ (which together with Standard conjecture D imply the previous mentioned two conjectures), proposed in previous works by Fei Hu and the author. The heuristic argument relies on the possibility of defining the self-composition $Fr_X^s$ in a good way, where $s$ is an arbitrary rational number (allowed to be negative), and similarly for another object related to the Frobenius morphism. 


\end{abstract}


The main point of this paper is to illustrate that Standard conjecture D (alone) is/should be enough to prove some useful and interesting open questions around Weil's Riemann hypothesis. This in turn may suggest some new ideas towards its analytic counterpart, the Riemann hypothesis.  

\subsection{A brief overview of Weil's Riemann hypothesis and some generalisations}

Let $X$ be a smooth projective variety defined over a finite field $\mathbb{F}_q$.  Let $Fr_X$ be the Frobenius morphism on $X$ (in fact, the discussion below is still valid if we replace $Fr_X$ by another polarised endomorphism on $X$).  We choose a Weil's cohomology theory $H^*(X)$, and for simplicity we assume moreover that the cohomological groups are vector spaces over a subfield of $\mathbb{C}$ (for example, \'etale cohomology), see \cite{Kleiman68}\cite{KM74}.  We recall the well known: 

{\bf Weil's Riemann hypothesis:} Let $\lambda$ be an eigenvalue of $Fr_X^*:H^k(X)\rightarrow H^k(X)$. Then $|\lambda |=q^{k/2}$. 

This conjecture inspired a lot of research in Algebraic Geometry since the 1960s. In particular, Grothendieck and Bombieri introduced the Standard conjectures (see \cite{Kleiman68}), modeled on a suggestion from \cite{Serre60}, whose validity will solve Weil's Riemann hypothesis. Surprisingly, Deligne \cite{Deligne74} solved Weil's Riemann hypothesis without solving the Standard conjectures. Interestingly, the proof of Deligne also yielded the proofs of several of the Standard conjectures. Nowadays, except the Standard conjecture C (which is a consequence of Deligne's proof, see \cite{KM74}) the other conjectures are largely unsolved: Lefschetz type standard conjecture (Standard conjecture A, B), Standard conjecture D and the Hodge Standard conjecture (a counterpart of the well known Hodge-Riemann relations on compact K\"ahler manifolds). If the Hodge standard conjecture holds, then Standard conjectuers A, B and D are equivalent, for many other relations between the Standard conjectures, the reader can consult \cite{Kleiman68}. In contrast to the situation over fields of zero characteristic, over finite fields the current status seems to indicate that the Hodge Standard conjecture is much harder than Standard conjectures A, B, D. In positive characteristic, besides easy examples (like products of projective spaces), Standard conjecture D is proven for Abelian varieties \cite{Clozel99} and certain self-products of K3 surfaces \cite{IIK25}, and the Hodge standard conjecture is known for Abelian varieties of dimension $4$ \cite{Ancona21} and certain self-products of K3 surfaces \cite{IIK25}. Moreover, Standard conjecture D is preserved, with certain appropriate conditions, under blowups and taking quotient spaces (see Lemma 2.6 in \cite{HTX2025}  for details),  while it is unclear how Hodge Standard conjecture will behave. For the reader's convenience, we recall: 

{\bf Standard Conjecture $D$.} (Grothendieck and Bombieri)  
For every algebraic cycle $Z$ on $X$ of codimension $k$, if $Z$ is numerically trivial then it is homologically trivial. 

In this paper we concentrate on several generalisations of Weil's Riemann hypothesis, which we will recall now. The readers are referred to \cite{Truong25}\cite{HT23}\cite{HT24}\cite{Xie24}\cite{HTX2025} for more detail. 

First, we recall the semisimplicity question by Tate: 

{\bf Semisimplicity Conjecture.} (cf. \cite[\S3, Conjecture (d)]{Tate65})
Let $f:X\rightarrow X$ be an $a$-polarized endomorphism, i.e., there exist an ample divisor $H_X$ on $X$ and $a\in \bZ_{>1}$ such that $f^*H_X \sim a H_X$.
Then for each $i$, the pullback map $f^*: H^i(X)\rightarrow H^i(X)$ is semisimple, i.e. diagonalisable.

Tate also proposed (which is an analog of the mentioned result in \cite{Serre60})  that the statement of Weil's Riemann hypothesis also holds more generally  for all polarised endomorphisms of $X$. This analog of Serre's result is now a theorem of Junyi Xie \cite{Xie24}, see also \cite{HTX2025}.  The above mentioned Semisimplicity conjecture is regarded as among the most important current open problems in Algebraic Geometry, providing a deeper understanding of Weil's cohomology theories in positive characteristics and other applications; and roughly speaking, it is an analog to the conjecture that all zeros of the Riemann xi function are simple. 

The author (inspired by work of Esnault and Srinivas  \cite{ES13}) has proposed yet more general conjectures. Let us first recall some backgrounds.  A dynamical correspondence on $X$ is an effective algebraic cycle $V$ of codimension $n$ in $X\times X$, so that the natural projections to the first and second factors of each irreducible component $W$ are surjective. A dynamical correspondence is finite (developed by Suslin and Voevodsky, see \cite{MVW06}) if these natural projections are finite morphisms on each irreducible component. For two dynamical correspondences $f$ and $g$ on $X$, we can compose them dynamically $f\Diamond g$, to obtain a new dynamical correspondence. We refer the reader to \cite{DS08, Truong20, HT24} for an overview on this dynamical correspondence, where also the relation to the usual composition of algebraic cycles in Algebraic Geometry is discussed. To illustrate, if $f$ and $g$ are dominant rational maps, then $f\Diamond g$ is the same as the usual composition of dominant rational maps. If $f$ is a dynamical correspondence and $t$ is a positive integer, we denote by $f^{\Diamond t}=f\Diamond \ldots \Diamond f$ (t times) the t-th iterate. As endomorphisms, dynamical correspondences also pullback on cohomology and numerical groups. 

For a finitely dimensional complex vector space $H$, we choose a norm $||.||_H$ (and for convenience will write $||.||$ if no confusion arises). Since any two norms on $H$ are equivalence, for what discussed next a specific choice of $||.||$ does not make any difference. 

Let $f:X\vdash X$ be a dynamical correspondence. \cite{Truong20} shows that the following limit  $\lim _{t\rightarrow\infty}||(f^{\Diamond t})^*|_{N^k(X)_{\mathbf{R}}}||^{1/t}$ exists, and hence we can denote it $\lambda _k(f)$ (which is a generalisation of the spectral radius of $f^*|_{N^k(X)}$ in case $f$ is an endomorphism).  An alternative approach, for rational maps $f$, is given in \cite{Dang20}. These are the counterparts of those on compact K\"ahler manifolds, established earlier in  \cite{DS05a, DS08}.  A sufficient condition to guarantee that the sequence $\lambda _j(f)$, $j=0,1,\ldots ,n$ is log-concave is that there is an infinite sequence $n_l$ such that the graphs of the iterates $f^{\Diamond n_l}$ are irreducible, see \cite[Theorem 1.1]{Truong20}. This is in particular the case if $f$ is a dominant rational map on $X$ (in this case, see also \cite{Dang20}). 

In \cite{Truong25}, a generalisation of $\chi _k(f)$ to dynamical correspondences is also given. For $k=0,\ldots ,n$, we define: 
$$\chi _k(f)=\limsup _{t\rightarrow \infty}||(f^{\Diamond t})|_{H^k(X)}||^{1/t}.$$
It is unknown whether $\chi _k(f)$ is always finite, or if the limsup in the definition can be replaced by lim. On the other hand, \cite{Truong25} shows that if Standard conjecture D holds, the natural extension of Entropy comparison conjecture (i.e. $\max _{j=0,1,\ldots ,2 \dim (X)}\chi _j(f) =\max _{k=0,1,\ldots , \dim (X)}\lambda _k(f)$) to all dynamical correspondences also holds, and as a consequence we then have $\chi _k(f)<\infty$ for all $k$. 

By results in \cite{Truong25}, the following conjecture is a generalisation of Weil's Riemann hypothesis (which corresponds to the special case $f=$ the Frobenius morphism): 

\textbf{Dynamical degree comparison conjecture}: Let $f$ be a dynamical correspondence on $X$. Then for all $j=0,\ldots ,\dim (X)$ we have $\chi _{2j}(f)=\lambda _j(f)$. 

In fact,  \cite{ES13} solves Dynamical degree comparison conjecture for the case $f$ is an automorphism of a surface $X$. They also mentioned therein a comment by Deligne that their result would follow from Standard conjectures. By the end of this note, we will provide a heuristic argument that Dynamical degree comparison conjecture, in any dimension and for any dynamical correspondence, should follow from Standard conjecture D (alone). Relevant to this, we note that it is proven in \cite{HTX2025} that if Standard conjecture D holds on $X\times X$,  then Dynamical degree comparison conjecture holds for all dominant rational maps of $X$. 

See also \cite{Truong25}\cite{HT24} for works on this and other related conjectures (see also \cite{Hu19,Hu24}\cite{Shuddhodan19}). In particular, the following conjecture from  \cite{HT24} implies both Dynamical degree comparison and Semisimplicity conjecture: 

\textbf{Norm comparison conjecture}: There is a constant $C>0$ such that for any dynamical correspondence $f$ of $X$, we have for any $0\leq j\leq \dim (X)$: 
\begin{eqnarray*}
||f^*|_{H^{2j}(X)} || \leq C ||f^*|_{N^j(X)}||, 
\end{eqnarray*}
and for $0\leq j\leq \dim (X)-1$
\begin{eqnarray*}
||f^*|_{H^{2j-1}(X)} || \leq C \sqrt{||f^*|_{N^j(X)}||\times ||f^*|_{N^{j+1}(X)}||}. 
\end{eqnarray*}

Moreover, Norm comparison conjecture allows to bound the growth  of the pullback of iterations of an endomorphism on cohomology groups in terms of that on algebraic cycles. More precise, assume that Norm comparison conjecture holds for an endomorphism $f$. Let $b_{coh}$ be the maximal size of Jordan blocks of all eigenvalues $\lambda $ of $f^*:H^{2j}(X)\rightarrow H^{2j}(X)$ with $|\lambda |=\chi _{2j}(f)$ means the spectral radius of a square matrix. We define similarly $b_{alg}$ for $f^*:N^j(X)\rightarrow N^j(X)$, where the eigenvalues $\lambda$ considered have $|\lambda |=\lambda  _{j}(f)$.  Since Dynamical degree comparison conjecture follows from Norm comparison conjecture, we obtain $\chi _{2j}(f)=\lambda _j(f)$. Then the growth of $||(f^n)^*|_{H^{2j}(X)}||=||(f^*)^n|_{H^{2j}(X)}||$ is comparable to $n^{b_{coh}}\lambda _j(f)^n$ and bounds that of  $||(f^n)^*|_{N^{j}(X)_{\mathbb{Q}}}||$ from above. On the other hand, the growth of  $||(f^n)^*|_{N^{j}(X)_{\mathbb{Q}}}||$ is comparable to $n^{b_{alg}}\lambda _j(f)^n$. Therefore, we obtain the following application (which is used in \cite{HT24} when $f$ is a polarised endomorphism, whence $b_{alg}=1$). 

\begin{proposition} Assume that Norm comparison conjecture holds for $X$. Let $f$ be an endomorphism of $X$, and define the numbers $b_{coh}$ and $b_{alg}$ as in the previous paragraph. Then $b_{coh}=b_{alg}$. 
\label{PropositionBoundJordanBlock}\end{proposition} 

The above Proposition can be extended to dominant rational maps which are algebraically stable, that is those for which $(f^{\Diamond n})^*=(f^*)^n$ for all positive integers $n$. 

In \cite{HT24}, it is shown that Norm comparison conjecture follows from Standard conjecture D and the following conjecture. For any $r\in \mathbb{Q}_{>0}$, we define $\gamma _r:H^*(X)\rightarrow H^*(X)$ to be such that on $H^j(X)$ it acts as multiplication by $r^j$. Since Standard conjecture C holds on finite fields (as a consequence of Deligne's proof of Weil's Riemann hypothesis), $\gamma _r$ is represented by an algebraic cycle $G_r$ with rational coefficients (but it may not be effective). 

\textbf{Conjecture $G_r$}: There is a constant $C>0$ so that for all for all dynamical correspondences $f$ of $X$:
\begin{equation}
||G_r\circ f|_{N^n(X\times X)_{\mathbf{R}}} || \leq C \deg (G_r\circ f),
\end{equation}
where $\deg (.)$ is with respect to a fixed ample divisor $H_{X\times X}$. 

Conjecture $G_r$ is known for Abelian varieties and some other special cases, see  \cite{HT24} for detail. More generally, Conjecture $G_r$ holds if the algebraic cycle $G_r$ is {\bf effective} (an analog to what discussed in the next paragraph, as well as in Section 0.3, is that the algebraic cycle representing $(Fr_X^*)^s$ is effective if $s$ is an integer). Indeed, if the cycle $G_r$ is effective, then we can prove the following stronger version of Conjecture $G_r$, where the concerned correspondences $f$ need only be effective (and not necessarily dynamical), which can have more applications: 

\textbf{Stronger Conjecture $G_r$}: There is a constant $C>0$ so that for all for all {\bf effective} correspondences $f$ of $X$:
\begin{equation}
||G_r\circ f|_{N^n(X\times X)_{\mathbf{R}}} || \leq C \deg (G_r\circ f),
\end{equation}
where $\deg (.)$ is with respect to a fixed ample divisor $H_{X\times X}$. 

 For $f=$\textbf{polarised endomorphisms}, it is shown in \cite{HT23} that Conjecture $G_r$ is a consequence of Standard conjectures.  However, it is unknown if Conjecture $G_r$ for general dynamical correspondences also follow from Standard conjectures. A main upshot of this paper is to heuristically show that nonetheless Norm comparison conjecture should be a consequence of Standard conjecture D, thus giving more support to its validity. We will also heuristically show that Stronger Conjecture $G_r$ should follow from Standard conjecture D. Also, in the last part of this paper, we will discuss a replacement of $\gamma _r$ in (Stronger) Conjecture $G_r$ by another cohomological operator which is more geometric in nature.  

Recently, Dynamical degree comparison conjecture is solved for endomorphisms $f$ in \cite{Xie24} (see also \cite{HTX2025} for another proof, as well as for some generalisations).  In particular, if $f$ is a polarised endomorphism, then the statement of Weil's Riemann hypothesis holds for $f$. One crucial idea in \cite{Xie24} is to estimate the norms of $(Fr_X^*)^s\circ (f^*)^t$ (as linear maps on the cohomology groups) for positive integers $t$ and for (possibly negative) integers $s$. It is \textbf{very important} here that $s$ is allowed to be negative. 

\subsection{Standard conjecture D implies Semisimplicity conjecture} Here we state and prove the first main result of this note. We first recall the following simple Linear Algebra fact, including a sketch of proof for the reader's convenience. If $A$ is a square matrix with complex coefficients, then we let $A^{\tau}$ denote the complex conjugate of $A$, $Tr(A)$ the trace of $A$, and $sp(A)$ the spectral radius of $A$. Recall that $Tr(A)$ is the sum of the diagonal entries of $A$, and $sp(A)=$ max $\{|\lambda|: \lambda $ is an eigenvalue of $A\}$. 

\begin{lemma} Let $A$ be a square matrix with complex coefficients. Then we have
$$sp(AA^{\tau })\leq (1+ sp(A))^2.$$
\label{LemmaBound2Norm}\end{lemma}
\begin{proof}[Sketch of proof] Let $A=UJU^{-1}$ be the Jordan normal form of $A$, where $U$ is an invertible matrix, and $J$ consists of Jordan blocks. Using the properties
$$sp(AA^{\tau })=\limsup _{n\rightarrow\infty}|Tr((AA^{\tau})^n)|^{1/n},$$
and $Tr(A_1A_2)=Tr(A_2A_1)$ for every two square matrices $A_1,A_2$ of the same size, we have $$sp(AA^{\tau })=\limsup _{n\rightarrow\infty}|Tr((AA^{\tau})^n)|^{1/n}=\limsup _{n\rightarrow\infty}|Tr((JJ^{\tau})^n)|^{1/n}=sp(JJ^{\tau}).$$

Therefore, we can assume, without loss of generality, that $A$ is itself a Jordan block $J$ with $\lambda$ on the diagonal. 

Let $N$ be the dimension of $J$. Then it is easy to check that $B=JJ^{\tau}$ is a Hermitian, positive semidefinite $N\times N$ matrix whose entries $b_{ij}$'s are as follows: 

$b_{i,j}=0$ if $|i-j|>1$, 

$b_{i,i}=1+|\lambda |^2$ for $i=1,\ldots, N-1$,

$b_{N,N}=|\lambda |^2$, 

$b_{i,i+1}=\overline{\lambda}$ for $i=1,\ldots ,N-1$, 

$b_{i-1,i}=\lambda$ for $i=2,\ldots ,N$. 

Since $B$ is Hermitian and positive semidefinite, it follows that $sp(B)=\sup _{v\in \mathbf{C}^N: ||v||=1}<Bv,v>$, which - from the above description for the entries of $B$ - can be bound from  above by $(1+|\lambda|)^2=(1+sp(J))^2$.   
\end{proof}

Now we show that Standard conjecture D implies Semisimplicity conjecture. This strengthens Theorem 1.1 in  \cite{HTX2025} which asserts that if Standard conjecture D holds on $X\times X$ and there is one polarised endomorphism $f$ on $X$ which is semisimple, then all polarised endomorphisms on $X$ are semisimple. 

\begin{theorem} Assume that Standard conjecture D holds  (for algebraic cycles of dimension $=\dim (X)$) on  $X\times X$. Then Semisimplicity conjecture holds. 

\label{TheoremStandardConjectureDSemisimplicityConjecture}\end{theorem}
\begin{proof}
We let $f$ be an $a_1$-polarised endomorphism on $X$ ($a_1>1$). Fix $k\in \{0,\ldots ,\dim (X)\}$. Then for all  integer $s$ (maybe negative) and positive integers $t$, by Lemma 2.5 in  \cite{HTX2025} (or use directly arguments in the earlier work  \cite{Truong25}), it follows that for  $A=Fr^*:H^k(X)\rightarrow H^k(X)$ we have 
$$sp((A^{-s})^{\tau}A^{-s})^{-1/2}||(f^*)^t|_{H^k(X)}||\leq C \max _{j=0,\ldots ,\dim (X)}q^{sj}a_1^{tj}.$$
If $s\leq 0$, note that $A^{-s}$ (and hence also $(A^{-s})^{\tau}$) has all eigenvalues of absolute value $q^{-sk/2}\geq 1$. Therefore  (here is where the delicate of this proof lies), by Lemma \ref{LemmaBound2Norm} applied to $(A^{-s})^{\tau}$, we have  $sp(A^{-s}(A^{-s})^{\tau})^{-1/2}\geq q^{sk/2}/2$. Hence, we obtain, for all integers $s\leq 0$ and all positive integers $t\geq 0$: 
$$q^{sk/2}||(f^*)^t|_{H^k(X)}||\leq C \max _{j=0,\ldots ,\dim (X)}q^{sj}a_1^{tj}.$$

We can combine this with the proof of part 2 of Theorem 1.1 in  \cite{HTX2025} (note that therein only the case $s\leq 0$ is needed, via the use of Dirichlet's approximation theorem) to conclude that all the Jordan blocks of $f^*:H^k(X)\rightarrow H^k(X)$ have size $1$, i.e. Semisimplicity conjecture holds. {For the reader's convenience, here we provide the detailed argument.}

We let $\nu $ be the largest Jordan block of $f$. We will show that $\nu =1$. Since $f$ is $a_1$-polarised, it follows (from the mentioned results in  \cite{Xie24}, see also  \cite{HTX2025}) that $||(f^*)^t|_{H^k(X)}||$ is of the same size as $a_1^{tk/2} t^{\nu -1}$. Then, by the previous inequality, we have (for integers $s\leq 0$ and $t\geq 0$)
\begin{equation}
q^{sk/2}a_1^{tk/2} t^{\nu -1}\leq C \max _{j=0,\ldots ,\dim (X)}q^{sj}a_1^{tj}.
\label{EquationSemisimplicityInequality}\end{equation}

Since both $q,a_1>1$, the number $\alpha$, defined as
$$\alpha =-\log (a_1)/\log (q),$$   
is a {\bf negative} real number. The inequality (\ref{EquationSemisimplicityInequality}) becomes  (for integers $s\leq 0$ and $t\geq 0$):  
\begin{equation}
e^{(s-t\alpha )(k \log q/2)} t^{\nu -1}\leq C \max _{j=0,\ldots ,\dim (X)}e^{(s-t\alpha )(j \log q)}.
\label{EquationSemisimplicityInequality1}\end{equation}

By Dirichlet's approximation theorem, there is a sequence $t_m\rightarrow +\infty$ and $s_m\in \mathbf{Z}$ such that $$|t_m\alpha -s_m|<1/t_m$$ 
for all $m$. Since $\alpha <0$, it follows easily that $s_m$ is negative, and hence can be used in the inequality (\ref{EquationSemisimplicityInequality1}). With this choice of $s_m$ and $t_m$, then both $e^{(s_m-t_m\alpha )(k \log q/2)}$ and $\max _{j=0,\ldots ,\dim (X)}e^{(s_m-t_m\alpha )(j \log q)}$ are bounded. From the inequality  (\ref{EquationSemisimplicityInequality1}), it follows that also the sequence $t_m^{\nu -1}$ is bounded, which implies $\nu=1$ as wanted. 

\end{proof}

For previous results on relations between semisimplicity, Standard conjectures (specially Standard conjecture D) and Tate's conjecture on algebraic cycles,  in different contexts and viewpoints, see for example  \cite{Moonen19}, \cite{Fu99}, \cite{Milne86-AJM}, and \cite{Jannsen92}. It is known that Tate's conjecture on algebraic cycles imply Standard conjecture D.  

Remark that by the results in  \cite{Xie24} (see also  \cite{HTX2025}), in the above proof we can replace the use of $Fr_X$ by any polarised endomorphism of $X$. 

In the next Subsection, we will provide a heuristic argument to show that Standard conjecture D should imply also Norm comparison conjecture (which as said, is stronger than Semisimplicity conjecture). 

\subsection{Standard conjecture D should imply Norm comparison and Dynamical degree comparison conjectures}

In \cite{HTX2025}, it is observed that the approach in \cite{Xie24} has several similarities to that in  \cite{HT24}. We list here a couple of the most important ones: 

1. An estimate similar to Lemma 4.9 in   \cite{HT24} is established. 

2. While $\gamma _r$ is multiplication by $r^j$ on $H^j(X)$, $(Fr_X^*)^s$ acts (roughly, it is precisely so if $Fr_X^*$ is semisimple) as multiplication by $r^j$ after taking \textbf{norms} of elements in $H^j(X)$, provided that $s$ is an \textbf{integer} and here $r=q^{s/2}$.

It is noteworthy to mention here that \cite{Xie24} uses the fact that $Fr_X$ \textbf{commutes} with all correspondences on $X$. However, this is not necessary, the same argument applies if we replace $Fr_X$ by any other polarised endomorphism $F$ on $X$ so that we know the statement of Weil's Riemann hypothesis holds for $F$. (Recall that now this latter assumption is true unconditionally, by \cite{Xie24}.) 

However, there is a \textbf{very important difference} in this similarity. That is, for $\gamma_r$ we can define for every positive rational numbers $r$. However, for point 2 above we can only rigorously define  $(Fr_X^*)^s$ for \textbf{integers} $s$, and hence  for positive numbers $r$ of such forms $r=q^{s/2}$.  This has the inconvenience that we are forced to look for sequences of $s_m$ and need to compose with iterations of $f^{\diamond t_m} $, and then take limits in the estimates. When taking limits, we obtain weaker estimates than needed for the above conjectures (Semisimplicity, Dynamical degree comparison, and Norm comparison).  

What if we can define $(Fr_X^*)^s$ for all  \textbf{real numbers} $s$ (possibly negative) with formal properties preserved? Then, together with Standard Conjecture D we  can resolve both Dynamical degree comparison and more generally Norm comparison conjectures. Before stating the result, let us recall some preliminary facts. 

First, a convenient norm on $H^k(X)$ is as follows. We choose $v_1,\ldots ,v_N$ be a basis for $H^k(X)$, and $w_1,\ldots ,w_N$ the  dual basis on $H^{2\dim (X)-k}$ (with respect to the cup product). Then for $v\in H^k(X)$, we define  $$||v||:=\sum _{j=1}^N |v\cdot w_j|,$$
which by Poincar\'e duality can be easily checked to be a norm. With this norm, we can define a norm for correspondences as follows: if $f$ is a correspondence (i.e. an algebraic cycle of dimension $=\dim (X)$ on $X\times X$, not necessary effective or dynamical), then 
$$||f^*|_{H^k(X)}||:=\sum _{i=1}^N ||f^*(v_i)||=\sum _{i,j=1}^N |f^*(v_i)\cdot w_j|.$$ 
The above norm is explicitly computed as follows. Let $\pi _1,\pi _2:X\times X\rightarrow X$ be the projections to the first and second factors.  Then $f^*(v_i)=(\pi _1)_*(f\cdot \pi _2^*(v_i))$. Therefore, 
 $$||f^*|_{H^k(X)}||=\sum _{i,j=1}^N |f\cdot \pi _2^*(v_i)\cdot \pi _1^*(w_j)|.$$ 

There is a natural involution $\tau :X\times X\rightarrow X\times X$ defined as follows: $\tau (x_1,x_2)=(x_2,x_1)$. For a correspondence $f$, we define its transpose $f_{\tau}$ by the formula: $f_{\tau}=\tau (f)$. With the above norms, we can check that for all $k$: 
 $$||f^*|_{H^k(X)}||=||f_{\tau}^*|_{H^{2\dim (X)-k}(X)}||.$$   
Similarly, for all $j$ we have $||f^*|_{N^j(X)}||=||f_{\tau}^*|_{N^{\dim (X)-j}(X)}||$. 

We fix an inner product on the complex vector space $H^*(X)$, such that different generalised eigenspaces of the Frobenius morphism are orthogonal. With this, for any matrix $A$ (or linear operator on $H^*{X}$), we can define its complex conjugate $A^{\tau}$. Here is the second main result/approach of this note: 

\begin{theorem}
Assume that for every real number $s$ (allowed to be both positive and negative), we can define $(Fr_X^*)^s$ on $H^*(X)$ in such a way that: i) the eigenvalues of $(Fr_X^*)^s$ on $H^j(X)$ have absolute values $q^{sj/2}$;   and  ii) the conclusion of the statement of Lemma 2.5 in  \cite{HTX2025} holds for $t=1$, that is: 
\begin{equation}
sp((A^{-s})^{\tau}A^{-s})^{-1/2}||f^*|_{H^k(X)}||\leq C \max _{j=0,\ldots ,\dim (X)}q^{sj}||f^*|_{N^j(X)_{\mathbf{Q}}}||,
\label{Equation100}\end{equation}
where $A=Fr^*:H^k(X)\rightarrow H^k(X)$.

Then Norm comparison conjecture holds on $X$. Consequently, both Dynamical degree comparison conjecture and Semisimplicity conjecture holds. 
\label{TheoremStandardConjectureNormComparisonConjecture}\end{theorem}
 \begin{proof}
As established in  \cite{HT24}, we only need to prove Norm comparison conjecture. To this end (see also the arguments in  \cite{HT24}, given that the norm on $N^j(X)$ is \textbf{additive} while the norm on $H^j(X)$ is \textbf{sub-additive} for a sum of dynamical correspondences), we can assume that $f$ is \textbf{irreducible}. With this we have the following \textbf{important} property: the sequence $a_j=||f^*|_{N^j(X)}||$, which we can take to be $f^*(H^j).H^{\dim (X)-j}$, is \textbf{log-concave} (i.e. $a_k^2\geq a_{k-1}a_{k+1}$ for all $1\leq k\leq \dim (X)-1$).  
 
 Under the assumption of the Theorem,  we  deduce that there is a constant $C>0$ so that for all dynamical correspondence $f$, integers $0\leq k\leq  \dim (X)$ and real numbers $r>0$ we have
 \begin{equation}
 r^{k}||f^*|_{H^k(X)}|| \leq C \max _{j=0,\ldots ,\dim (X)} r^{2j}||f^*|_{N^j(X)_{\mathbb{R}}}||. 
\label{EquationNormBound}\end{equation}
Indeed, if $r\leq 1$, then choose the real number $s=\log _q r\leq 0$, the statement of Lemma 2.5 in  \cite{HTX2025} combined with the proof of Theorem \ref{TheoremStandardConjectureDSemisimplicityConjecture} establish (\ref{EquationNormBound}). For the case $r>1$, with $r_1=1/r <1$ and what already proven, applied to the transpose $f_{\tau}$ (which is also a dynamical correspondence and is irreducible), we obtain: 
  \begin{eqnarray*}
 r_1^{2\dim (X)-k}||f_{\tau}^*|_{H^{2\dim(X)-k}(X)}|| \leq C \max _{j=0,\ldots ,\dim (X)} r_1^{2\dim (X)-2j}||f_{\tau}^*|_{N^{\dim (X)-j}(X)_{\mathbb{Q}}}||. 
\end{eqnarray*}
Dividing $r_1^{2\dim (X)}$ from both sides of the above inequality, using $r_1=1/r$ and the relations between the norms of $f^*$ and $f_{\tau}^*$, we obtain  (\ref{EquationNormBound}) for $r>1$.  

If the constant $C$ in  (\ref{EquationNormBound}) is smaller than or equal to $1$, then Lemma 4.9 in \cite{HT24} proves the  Norm comparison conjecture under the assumptions of the theorem.  For the case $C>1$ we can just redefine the norm on $N^j(X)_{\mathbb{Q}}$ to be  $C||f^*|_{N^j(X)_{\mathbb{Q}}}||$ and reduce to the case $C=1$. 
  \end{proof}

 \textbf{Remarks:} 
 
 1. Condition i) in the theorem should be automatic, from what we know about the eigenvalues of the Frobenius morphism by Weil's Riemann hypothesis.   
 
 2. In the proof of Lemma 2.5 in  \cite{HTX2025}, it is essential that for $s\in \mathbb{Z}$ and a dynamical correspondence $f$, the correspondence $(Fr^*)^s\circ f^*$ is effective. We expect that when we can define $(Fr^*)^s$ for real numbers $s$ then this needed property is preserved. Another assumption which is needed in Lemma 2.5 in  \cite{HTX2025} is that Standard conjecture D holds on $X\times X$. 
 
 On the other hand, it is plausible that we may be able to prove both conditions i) and ii) in Theorem \ref{TheoremStandardConjectureNormComparisonConjecture} without needing Standard conjecture D, see also the discussions below for some supporting evidences. Then, we also obtain Semisimplicity conjecture.  
 
 3. By the results on absolute values of eigenvalues of polarised endomorphisms in \cite{Xie24} (see also  \cite{HTX2025} ),  we can replace $Fr_X$ by any other polarised endomorphism.  
 
 4. Even for varieties where all Standard conjectures hold, it is unclear if Conjecture $G_r$  (for general dynamical correspondences, not just polarised endomorphism) also holds. In this case, the approach through defining  $(Fr^*)^s$ for all real numbers $s$ seems to be more achievable. 
 
 5. From the above points, here is a more explicit description of a new approach towards resolving the Norm comparison conjecture and hence the Dynamical degree comparison conjecture: 
 
 i) Establish Standard conjecture D  (for algebraic cycles of dimension $=\dim (X)$)  on  $X\times X$. 
 
 ii) For a polarised endomorphism $F$ on $X$ (for example $F=Fr_X$), define $(F^*)^s$ for all \textbf{real numbers} $s$.  Establishing that for all dynamical correspondences $f$, then $(F^*)^s\circ f^*$ is represented by an effective algebraic cycle,  or alternatively establishing an estimate similar to that in Lemma 2.5 in \cite{HTX2025}.   
 
\textbf{A relation between Standard conjecture D and Norm comparison and Dynamical degree comparison conjectures}: As a consequence of what discussed in points 4 and 5 above, we see that  the validity of Standard conjecture D should imply the validity of Norm comparison conjecture and Dynamical degree comparison conjecture. In fact, it remains to define in a natural way $(F^*)^s$ for all real numbers $s$, which should automatically satisfy condition i) in the theorem (since for $s=1$ or more generally an arbitrary integer then these conditions are valid). Condition ii) is, as explained, follows from that $(F^*)^s\circ f^*$ should be effective for all dynamical correspondences $f$. Again, this latter requirement holds for all integers $s$, with involved constants \textbf{independent} of $s$.  
 
\textbf{Some discussions pertaining a possible definition of $(Fr_X^*)^s$ for real numbers $s$}: We assume that $Fr_X$ is semisimple, i.e. $Fr_X^*$ is a a diagonalisable matrix (for example, if Standard conjecture holds on $X\times X$). Hence, we can write $Fr_X^*=UDU^{-1}$ where $D$ is an invertibly diagonal matrix. If $s$ is an integer, then $(Fr_X^*)^s$ is $UD^sU^{-1}$. Therefore, we will try to define $(Fr_X^*)^s=UD^sU^{-1}$, provided we know how to define $D^s$. Recall that we need to define $(Fr_X^*)^s$ in such a good way that the composition $(Fr_X^*)^s\circ f^*$ is close to be effective, so the estimates in Lemma 2.5 in \cite{HTX2025} are valid. We also note that we need only to do this for $s\in \mathbf{Q}$. To this end, it remains to define for the case $s=1/N$ where $N$ is a positive number. See also the next Subsection for a clearer idea when we assume all Standard conjectures.   

\subsection{Standard conjecture D should imply (Stronger) Conjecture $G_r$} We now provide a clearer and improved idea to that in the previous subsection, which provides support to Stronger Conjecture $G_r$ (which is stronger than both Norm comparison conjecture and Dynamical degree comparison conjecture). This is done through a direct relation between the cohomological operator $\gamma _r$ and the Frobenius morphism, see a relevant (less direct) discussion in the last part of \cite{HTX2025}. 

We revisit the previous section. The matrix $B=AA^{\tau}$ is positive definite Hermitian, therefore it has a {\bf unique} positive definite Hermitian $k$-th root (see e.g. Theorem 7.2.6 in \cite{HornJohnson2012}), for all positive integers $k$: We simply decompose $B=UDU^{-1}$, for $D$ a diagonal matrix with positive real numbers on the diagonal and $U$ a unitary matrix, and define $B^{1/k}=UD^{1/k}U^{-1}$. From this, it is ready to show that $B^{s}$ is uniquely defined as a positive definite Hermitian matrix, for all $s\in \mathbf{Q}$ (negative numbers are allowed). We can replace the use of $Fr^s$ in the discussion in the previous Subsection by $B^{s/2}$, the advantage is now $B^{s/2}$ is legitimately defined.  Though, there is still a question that which inner product on $H^*(X)$ should we choose. Here a good choice (to be clear why in the discussion in the next paragraphs) seems to be: We choose an inner product product on $H^*(X)$ such that the {\bf different generalised eigenspaces} of the Frobenius morphisms are {\bf orthogonal}.  We can ask the following question: 

{\bf Question:} Is $B^s$ independent of the choice of such an inner product on $H^*(X)$?

When Standard conjecture D holds, we will see next that the answer to this Question is affirmative. Moreover, it suggests us to ask the following version of Conjecture $G_r$:

{\bf Conjecture $B^s$}: 1) For all $s$, the cohomological operator $B^s$ can be represented by an algebraic cycle. 2) There is a constant $C>0$ such that for all dynamical correspondences $f$ we have
\begin{eqnarray*}
||B^s\circ f|_{N^n(X\times X)_{\mathbf{R}}} || \leq C  \sum _{k=0,\ldots ,\dim(X)} q^{ks}\deg _k(\circ f). 
\end{eqnarray*}

We can also obtain a stronger conjecture by requiring that the above inequality holds for all effective algebraic cycles $f$, and not just dynamical ones. 

The  validity of Conjecture $B^s$ plus Standard conjecture D again will imply both Norm comparison and Dynamical degree comparison conjectures, by replacing $G_r$ in the proofs in  \cite{HT24} by $B^s$ for $s=q^r$, and combining with the trick used in the proof of Theorem \ref{TheoremStandardConjectureNormComparisonConjecture} (consider separately the cases $s<0$ and $s>0$). In the above discussion, we can also use another polarised endomorphism in place of $Fr_X$, and the whole discussion goes through. However, since the Frobenius morphism is a natural morphism on $X$, it seems that we have more chance to proceed when using it instead of another polarised endomorphism. On the other hand, if Standard conjecture D holds, then essentially there is no difference between using any of the polarised endomorphism on $X$, see the next discussion.  

From now on, assume again that Standard conjecture D holds on $X\times X$. Then, by Theorem \ref{TheoremStandardConjectureDSemisimplicityConjecture}, $Fr_X^*:H^*(X)\rightarrow H^*(X)$ is diagonalisable. Therefore, we can choose a basis for $H^*(X)$, as a complex vector space, such that $A=Fr_X^*$ is a diagonal matrix. With this basis, comes a natural inner product $<.,.>$,  and we can then define the complex conjugate $A^{\tau}$.  Moreover, on $H^j(X)$, then $B^s$ acts as multiplication by $q^{js}$. Hence, clearly, with $r=q^s$ for $s\in \mathbf{Q}$, $B^s$ is the same as the cohomological operator $\gamma _r$ considered in (Stronger) Conjecture $G_r$. Therefore, $B^s$ is {\bf unambiguously} defined (and even $A^{\tau}$ itself is unambiguously determined, see later).  Therefore, $B^s$ provides a {\bf construction with geometric origin} for the (apparently) purely cohomology operator $\gamma _r$, which is moreover compatible with a complex structure on $H^*(X)$. This gives more support to the validity of Conjecture $G_r$, and below we will present some more explicit ideas on this. If we replace $Fr_X$  by another $a$-polarised endomorphism $F$, then the results in  \cite{Xie24} (see also \cite{HTX2025}) shows that the correspondingly constructed $B$ satisfies $B^s=\gamma _{a^s}$, and hence essentially there is no difference.  

Now if $B^{s}=\gamma _{q^s}$ is   the class of an {\bf effective} algebraic cycle, for all $s\in \mathbf{Q}$, then it follows from what mentioned before the statement of Stronger Conjecture $G_r$ that Stronger Conjecture $G_r$ is valid.  We note that the fact that $B^s$ is an algebraic cycle (i.e. representable by an algebraic cycle, with real coefficients) is clear, since it is the same as $\gamma _r$, and we saw that $\gamma _r$ is an algebraic cycle. Therefore,  only the effectiveness of $B^{s}$  is needed to be shown. We note that if this is true, then it seems a significant  fact.  

We note also that $A^{\tau}$ is itself an {\bf algebraic cycle} (with complex coefficients), under the assumption that Standard conjecture D holds on $X\times X$. This is based on the observation that, while the equality $<Av,w>=<v,A^{\tau}w>$ depends on the specific choice of a basis for  $H^*(X)$ by eigenvectors, $A^{\tau}$ - as a linear operator on $H^*(X)$ - is i{\bf independent} of this choice. Indeed, if we let $\Lambda $ be the set of eigenvalues of $A$, and for each $\lambda \in \Lambda$ denote by $E_{\lambda}\subset H^*(X)$ the corresponding eigenspace, then 
$$A^{\tau}=\sum _{\lambda \in \Lambda}\overline{\lambda} pr_{\lambda},$$ 
 where $pr_{\lambda}:H^*(X)\rightarrow E_{\lambda}$ is the linear projection. These linear projections $pr_{\lambda}$ are {\bf unambigously} defined (i.e. independent of the choice of a basis as above), by Hamilton-Cayley lemma, see e.g.  \cite{KM74}. Moreover, they are represented by linear combinations of graphs of iterates of the Frobenius morphism (with complex coefficients), and hence are algebraic cycles with complex coefficients. 
  
Even if $\gamma _r$ turns out to be not effective, we recall that $B$ is the composition of an effective algebraic cycle and its complex conjugate (which is itself an algebraic cycle with complex coefficients), $B$ is itself an algebraic cycle (with rational coefficients), and  $B^{s}$ is uniquely defined from $B$.  This may  provide some special properties which allow us to prove (Stronger) Conjecture $G_r$.  It is also interesting to find new examples of varieties $X$ for which $\gamma _r$ is effective, besides known cases as Abelian varieties (see \cite{HT24}). 

{\bf Remark}: We recall that Conjecture $G_r$ combined with Standard conjecture D  imply both Norm comparison and Dynamical degree comparison conjectures, by results in \cite{HT24}. Therefore, here we have another route to see that Standard conjecture D should imply both  Norm comparison and Dynamical degree comparison conjectures. 

 \subsection{Acknowledgements} The first two Subjections of this note were done while  the author was visiting Department of Mathematics at University of Tokyo and Center for Mathematical Science and Artificial intelligence at  Chubu University, in August 2025, and the last Subsection of this note is done while the author is visiting Department of Mathematics and Computer Science at University of Ferrara. The author would like to thank Keiji Oguiso,  Takayuki Watanabe and Cinzia Bisi, and these institutions for the hospitality and inspiring enviornments.  He also would like to thank the Department of Mathematics (University of Oslo),  University of Ferrara and INdAM  for financial support.

\bibliographystyle{amsalpha}
\bibliography{mybib}

\end{document}